\documentclass[a4paper,12pt]{amsart}

\usepackage{amsmath,amssymb}
\usepackage{algorithm2e}
\usepackage{graphicx}
\usepackage{booktabs}
\usepackage{multirow}
\usepackage{url}

\newcommand{\Ref}[1]{\mbox{(\ref{#1})}}
\DeclareMathAlphabet\mathbold{OML}{cmm}{b}{it}
\numberwithin{equation}{section}

\renewcommand{\div}{\operatorname{div}}
\newcommand{\Range}{\operatorname{Range}}

\def\b1{{\mathbf 1}}

\def\bu{{\mathbf u}}
\def\bw{{\mathbf w}}

\def\be{{\mathbf e}}

\def\bx{{\mathbf x}}

\newcommand{\R}{{\mathcal R}}
\newcommand{\A}{{\mathcal A}}

\newcommand{\M}{{\mathcal M}}

\def\b1{{\mathbf 1}}
\def\b0{{\mathbf 0}}

\def\bx{{\mathbf x}}
\def\bw{{\mathbf w}}
\def\be{{\mathbf e}}
\def\bu{{\mathbf u}}

\renewcommand{\Range}{{\mathcal Range}}
\renewcommand{\A}{{\mathcal A}}
\renewcommand{\M}{{\mathcal M}}
\DeclareMathOperator\blockdiag{blockdiag}
\DeclareMathOperator\diag{diag}
\renewcommand{\div}{\text{div}}

\begin{document}

\title[Improving solve time of Adaptive AMG]
{Improving Solve Time of Aggregation-based Adaptive AMG}
\author{Pasqua D'Ambra}
\address{Institute for Applied Computing ``Mauro Picone'', National Research Council of Italy, Naples, Italy}
\email{pasqua.dambra@cnr.it}

\author{Panayot S. Vassilevski}
\address{Center for Applied Scientific Computing, Lawrence Livermore National Laboratory, CA, USA\newline
$\text{\hspace{0.5cm}}$
Fariborz Maseeh Department of Mathematics and Statistics, Portland State University, Portland, OR, USA}
\email{panayot@llnl.gov, panayot@pdx.edu}

\thanks{This work of the first author was supported in part by EC under the Horizon 2020 Project
{\em  Energy Oriented Center of Excellence: toward exascale for energy – EoCoE II}, Project ID: 824158 and by the INdAM-GNCS Project on \emph{Innovative and Parallel Techniques for Large Linear and Nonlinear Systems, Function and Matrix Equations with Applications}.
This work of the second author was performed under the auspices of the
  U.S. Department of Energy by Lawrence Livermore National Laboratory
  under Contract DE-AC52-07NA27344 and was also supported in part by NSF under grant DMS-1619640.}

\begin{abstract}
This paper proposes improving the solve time of a bootstrap AMG designed previously by the authors.
This is achieved by incorporating the information, set of {\em algebraically smooth} vectors, generated by the bootstrap algorithm, in a single hierarchy by using sufficiently large aggregates, and these aggregates are compositions of aggregates already built throughout the bootstrap algorithm.
The modified AMG method has good convergence properties and shows significant reduction in  both, memory and solve time. These savings with respect to the original bootstrap AMG are illustrated on some difficult (for standard AMG) linear systems arising from discretization of scalar and vector function elliptic partial differential equations (PDEs) in both 2d and 3d.
\end{abstract}


\maketitle
\section{Introduction}\label{section: introduction}

Several black-box Algebraic MultiGrid (AMG) solvers have been introduced in the near past (\cite{BFMMMR-06,BBKL-11,DV-13}) which avoid any a priori assumptions on the type of matrices or their origin.
They exploit hierarchies of coarse vector spaces and respective matrices computed variationally, in the symmetric positive definite (s.p.d.) case, from the given, original, fine-level sparse matrix. All these methods share the common feature, that a current method (initially, a single-level one, such as Gauss-Seidel) is tested on the homogeneous problem $A \bx = \b0$ with non-zero initial guess and if convergence is deemed unacceptable, the most recent iterate (referred to as {\em algebraically smooth} error component) is incorporated in the existing hierarchy of vector spaces. Then the test is performed again with the new, adapted solver.
The way the adaptation is done differs in the various proposed {\em adaptive}, also referred to as {\em bootstrap}, AMG methods.
The building (setup) phase of these algorithms can be fairly expensive, which is to be expected, but can be to a large extent amortized if the solver is used multiple times, e.g., for different right-hand sides or applied to near-by matrices. Therefore, to make these methods of practical interest it is important to have their solve phase reasonably fast, as discussed in~\cite{BBKL-15}.
The latter is the main objective of the present paper.

The method  proposed by the authors in~\cite{DV-13} and extended later in~\cite{DFV-18}, consists of a sequence of aggregation-based AMG cycles, denoted by $B_r$, $r=1,\dots$, where each consecutive $B_r$ is constructed on the basis of the product iteration matrix:
\begin{equation}
\label{comp_AMG}
E_{r-1} = (I-B_{r-1} A) \ldots (I-B_1A),
\end{equation}
which uses the previously constructed cycles $B_s$, $s < r$. Here $A$ is the given $n \times n$ sparse s.p.d. matrix.
This product (or rather its symmetrized version) is applied to a random vector several times to estimate its convergence radius. If the latter is deemed unacceptable, i.e., staying above certain desired convergence factor, the current iterate is used to create a hierarchy of coarse vectors spaces and interpolation matrices that relate any two consecutive levels of the hierarchy. The coarsening is based on the current-level algebraically smooth vector. On the first (fine) level this is the vector created by the multiplicative error operator $E_{r-1}$.
On any given level, the constructed algebraically smooth vector $\bw$  is used to create aggregates, by several steps of pairwise aggregation.
Each step of pairwise aggregation uses the vector to define edge weights of the matrix graph which are then used in a specialized weighted matching algorithm.
Then, a piecewise-constant interpolation matrix $P$ is formed, by projecting the current smooth vector on the aggregates, and the coarse matrix $A_c = P^TAP$ is computed.
On each consecutive coarse level, the restriction of the previous (fine) level algebraically smooth vector $\bw$, i.e., $\bw_c = P^T \bw$, defines the next-level algebraically smooth vector  $\bw:= \bw_c$.
Applying this process recursively, until a reasonable small coarse size is reached, the hierarchy of vector spaces is constructed. They, together with standard (e.g., Gauss-Seidel) smoothers define the (unsmoothed) aggregation based AMG operator $B_r$.

In summary, each $B_r$ is constructed on the basis of an algebraically smooth vector $\bw_r$ which is in turn constructed on the basis of the previous operators $B_s$, $s < r$.
There is a clear analogy with Krylov-type methods; indeed, in Krylov methods, some search vectors are generated, where at every new iteration step, a new search vector is constructed based on the previous search vectors and the (preconditioned) residual. In the above described adaptive AMG, at each step, a new algebraically smooth vector is constructed based on the previous ones (via their respective operators $B_r$) and the role of preconditioning is played by the construction of the hierarchy leading to $B_r$.
Both the Krylov methods and the adaptive AMG will converge in at most $n$ steps, or equivalently, will generate at most $n$ vectors (in exact arithmetic), and in practice we expect this to take much smaller number of steps.

We note that the sequence $\{B_r\}$ requires storing of substantial amount of information, since each $B_r$ comes with a hierarchy of matrices $\{A^{k}_r\}^{nl}_{k=0}$ and respective interpolation matrices $\{P^{k}_r\}^{nl-1}_{k=0}$, where $nl$ generally depends on $r$.
The goal of the present paper is to reduce this large memory requirement, which is a key issue for scalable AMG on modern parallel architectures having millions of cores and decreased memory per core (see e.g., ~\cite{BGSY-11,DF-16}). Since the vectors $\{\bw_r\}$, via their respective $\{B_r\}$, capture (by construction) all components of the error, the idea is that
a single hierarchy, that simultaneously coarsens these vectors, may as well do a similar job.
Savings in memory is also expected to come from the fact that although $\{\bw_r\}$ are linearly independent (by construction), their local versions, i.e., their restrictions on aggregates, may be
linearly dependent and this will let us use few of them locally. This however comes with additional cost, for example, if we use SVD to extract such linearly independent subset locally.
Furthermore, the resulting interpolation matrix will be denser which reflects also the sparsity of the coarse matrices.
These considerations are addressed in the remainder of the present paper, and our conclusions are drawn based on running several tests on difficult (for classical AMG) matrices.
For these test problems the \emph{single-hierarchy multiple-vector aggregation-based AMG} resulting from our modification of the composite adaptive AMG shows good convergence and superior performance; it requires less memory and has faster solution time.
We should stress though, that even with the achieved improvement over the original bootstrap AMG, the resulting method can still be fairly expensive compared to the more traditional AMG ones.
However, the objective here is to design a genuinely {\em black-box} algebraic solver without any a priory assumptions and avoiding the use of any additional information that may be otherwise  available for a particular problem at hand.

The remainder of the paper is organized as follows.
In Section~\ref{section: BAMG algorithm}, we briefly summarize the algorithm for generating algebraically smooth vectors.
Section~\ref{section: review of compatible weighted matching coarsening} describes the way we generate aggregates by weighted matching where the weights come from a given algebraically smooth vector.
Section~\ref{section: multiple vectors interpolation} provides the details on the construction of the multiple vector interpolation matrices in aggregation AMG.
The assessment of the quality of our modified multiple-vector aggregation-based AMG is found in Section~\ref{section: numerical results}.
We close with some conclusions in Section~\ref{section: conclusions}.

\section{Generating algebraically smooth vectors using composite AMG} \label{section: BAMG algorithm}

In~\cite{DV-13,DFV-18} we proposed a new adaptive AMG method employing a bootstrap process
which has the final objective to setup a composite AMG as in~\Ref{comp_AMG} with a prescribed convergence factor.
Each $B_r$ is an AMG operator built by an aggregation procedure, named {\em coarsening based on compatible weighted matching}, which uses the most recent sample of smooth vectors dynamically generated by testing the last available composite AMG on the homogeneous linear system $A \bx = \b0$, starting from a random initial guess. Therefore, our bootstrap procedure (see Algorithm~\ref{alg-boot}), at each stage $r$, sets up a new AMG operator, i.e., a new hierarchy of aggregates and interpolation operators, and a new sample of smooth vector. The bootstrap can be stopped either if the composite AMG in~\Ref{comp_AMG} has reached the desired convergence rate or when a given number of hierarchies is built.
\begin{algorithm}[h]
\SetAlgoLined
\KwData{$A$: matrix, ${\bw_0}$: arbitrary vector, $\rho_{des}$: desired conv. rate, $maxstage$: max number of AMG operators,
$\nu$: number of iterations for testing phase}
\KwResult{AMG hierarchies to define $B_r$, smooth vectors $\bw_r$, $r=1, \ldots$}
{\bf Initialize:} $r=1$, $\rho=1.0$\;
\While{$\rho \geq \rho_{des}$ and $r \leq maxstage$}
{
\vspace{0.2cm}
{\bf Building Phase:} build new AMG operator $B_r$\;
\hspace*{1cm} apply {\em coarsening based on compatible weighted matching} to $A$ and ${\bw_{r-1}}$\;
{\bf Testing Phase:} {apply composite AMG and compute new smooth vector}\;
\hspace*{1cm} let $\bx^0$ a random vector\;
\hspace*{1cm} apply $\bx^j= \prod_{r} (I-B_r^{-1}A) \bx^{j-1},  \ j=1, \ldots, \nu$ (or a symmetrized version)\;
\hspace*{1cm} estimate convergence factor $\rho$;\\
\hspace*{1cm} $\bw_r=\bx^\nu/\|\bx^\nu\|_A$\;
}
\caption{Bootstrap Algorithm}
\label{alg-boot}
\end{algorithm}
The above procedure differs from the bootstrap AMG algorithm introduced in~\cite{BBKL-11} in the way it incorporates new information on smooth vectors in the final AMG. Indeed, our adaptivity approach re-applies the same aggregation algorithm for computing multiple
hierarchies of coarse spaces, each one of which includes one new smooth vector.
The application of a multiplicative composition of the computed hierarchies allows to obtain an efficient AMG method.
The approach in~\cite{BBKL-11}, at each bootstrap iteration, modifies an already computed AMG, characterized by a given coarse space and a corresponding prolongator, in order to include global basis of new smooth vectors in the previous coarse space by re-computing interpolation on all levels.

As observed in~\cite{BBKL-15}, due to the cost of the setup algorithm, a bootstrap AMG and more generally an adaptive AMG, is well suited for problems where standard AMG methods, using a priori characterization of smooth error components, lose their efficiency. Furthermore, the cost of the setup can be amortized when the same system has to be solved for multiple right-hand sides (r.h.s.), as in time-dependent problems, or when the method has to be applied to near-by matrices, as in Multilevel Monte Carlo simulations.

In the following we give a brief outline of the {\em coarsening based on compatible weighted matching}, which is our basic algorithm to build hierarchies of aggregates and prolongators to define each new operator $B_r$ in Algorithm~\ref{alg-boot}.

\section{Review of compatible weighted matching coarsening}\label{section: review of compatible weighted matching coarsening}

The basic kernel of Algorithm~\ref{alg-boot} is a completely automatic and general aggregation procedure for AMG methods which finds its starting point in the concept of {\em compatible relaxation} introduced in~\cite{B-00} as {\em a smoother that keeps the coarse variable invariant}. The above concept is widely used to compute measures of the suitability of a coarse space for a given problem~\cite{B-00,L-04,FV-04,BF-10,BBKL-15}.
In our method, we define a compatible relaxation by a
$D-$orthogonal decomposition of the original, fine-level, vector space $\R^n=\Range(P_c) \oplus^\perp \Range(P_f)$, where $D$ is a suitable s.p.d. matrix, by exploiting weighted matching in the matrix graph.  Let $G_C=(V,E,C)$ be the weighted undirected graph associated to
the symmetric matrix $A=(a_{ij})_{i,j=1, \ldots, n}$, where $V=\{1,2,\ldots,\;n\}$ is the index set, $E$ is the edge set and $C=(c_{ij})_{i,j=1, \ldots, n}$ is a matrix of non-negative edge weights.
A \emph{matching} in $G_C$ is a subset of edges $\M \subseteq E$ such that no two edges share a vertex, while a {\em maximum product matching} is a matching $\M$ that maximizes the product of the weights $c_{ij}$ of all edges $(i,j)\in \M$, i.e., $\arg \max_{\M} \prod_{ (i,j) \in \M} c_{ij}$.
Maximum product matchings are successfully used in sparse linear algebra to move large matrix
entries onto the main matrix diagonal (see~\cite{DK-01,HS-06,HS-15}).

Let  $\M=\{ e_1, \ldots, e_{n_p} \}$ be a matching of the graph $G_C$, with $n_p$ the number of index pairs, and let $\bw = (w_i)_{i=1, \ldots, n}$ be a given (smooth) vector; for each edge $e=(i,j)$, we can define the following two local vectors:
\begin{equation}\label{vector w_e}
 \bw_e =
   \frac{1}{\sqrt{w^2_i + w^2_j}}
\left [
\begin{array}{c}
w_i \\
w_j
\end{array} \right ] \quad \text{and} \quad
\bw^{\perp}_e =
\frac{1}{\sqrt{w^2_j/a_{ii} + w^2_i/a_{jj}}}
\left [
\begin{array}{c}
-w_j/a_{ii}\\
w_i/a_{jj}
\end{array} \right ],
\end{equation}
where $D_e= \left [
\begin{array}{cc}
a_{ii} & 0\\
0 & a_{jj}
\end{array} \right ]$ is the diagonal of the restriction of $A$ to the edge $e$.
Based on the above vectors, two {\em prolongators} can be constructed:
\begin{equation}
P_c = \left (
\begin{array}{cc}
\tilde{P}_c & 0\\
0 & W
\end{array}
\right ) \in \R^{n \times n_c}, \quad \text{and} \quad
P_f=
\left (
\begin{array}{c}
\tilde{P}_f \\
0
\end{array}
\right ) \in \R^{n \times n_p},
\label{prol}
\end{equation}
where:
\[
\tilde{P}_c=\blockdiag( {\bw}_{e_1}, \ldots, {\bw}_{e_{np}} ),
\quad \quad \tilde{P}_f=\blockdiag({\bw}^\perp_{e_1}, \ldots, {\bw}^\perp_{e_{np}}).
\]
$W=\diag(w_k/|w_k|), \; k=1, \ldots, n_s$, is related to possible unmatched nodes in the case $\M$ is
not a perfect matching for $G_C$ and $n_c=n_p+n_s$.

The matrix $P_c$ represents a piecewise-constant interpolation operator whose range includes the original (smooth) vector $\bw$; furthermore, by construction $(P_c)^TD P_f=0$, i.e., $\Range(P_c)$ and $\Range(P_f)$ are orthogonal with respect to the $D-$inner product on $\R^n$, with $D=diag(A)$. Exploiting the above decomposition, the matrix $A$ admits the following two-by-two block form:
\begin{equation}
[P_c,P_f]^T A [P_c, P_f]=
\left (
\begin{array}{ll}
P_c^TAP_c & P_c^TAP_f\\
P_f^TAP_c & P_f^TAP_f
\end{array}
\right ) =
\left (
\begin{array}{ll}
A_c & A_{cf} \\
A_{fc} & A_f
\end{array}
\right ).
\label{matdecomp}
\end{equation}
In the above setting, given a smoother $M$, the relaxation scheme defined by the following error propagation matrix:
\begin{equation}
E_f \equiv (I-P_f(P_f^TMP_f)^{-1}P_f^TA)
\label{comp-matrix}
\end{equation}
is a compatible relaxation for the two-level method whose coarse variables are defined by the prolongator $P_c$, i.e., $\bx_c=P_c^T \bx$. We observe that the sparsity pattern of the above prolongator is completely defined by the pairwise aggregation stemming from the matching $\M$
(see Algorithm 1 in~\cite{DFV-18}), while its values depend on the (smooth) vector $\bw$.
Relaxation defined by the matrix~\Ref{comp-matrix} is equivalent to an iteration of the form: $\be_{\ell+1}=(I-M_f^{-1}A_f) \be_\ell, \; \ell=1, \ldots,$ with $M_f = P^T_fM P_f$.
Therefore, when matrix $A_f$ is well conditioned or diagonally dominant, i.e., a Richardson-type relaxation method on $A_f$ is fast convergent, the coarse variables defined by the prolongator $P_c$ can be considered a suitable coarse set for an efficient two-level method~\cite{FV-04,XZ-17}. To this end we choose $\M$ so that the product of the diagonal entries of $A_f$ is as large as possible.
Note that $A_f$ has the same sparsity pattern of $A$ with entries depending on the entries of $A$ and of the vector $\bw$. It can be well represented by the undirected weighted graph $G_C=(G,C)$, where $G$ is the unweighted undirected graph of $A$ and $C$ is a suitable matrix of edge weights computable with linear complexity by simple algebraic operations for relations in~\ref{vector w_e}-\ref{matdecomp}.

In the above setting, we define the {\em coarsening based on compatible weighted matching}, described in details in~\cite{DFV-18}, that, starting from a maximum product matching of the weighted graph $G_C$, defines a general aggregation-based AMG which does not use any a priori information on the system matrix.
We observe that aggressive coarsening, with aggregates merging multiple pairs and having almost arbitrary large size of the type $n_c=2^s$ for a given $s$, can be obtained by combining multiple steps of the basic pairwise aggregation, i.e., by computing the product of ($s$) consecutive pairwise prolongators.  On the other hand, as expected, results discussed in~\cite{DV-16} show that the straightforward product of piecewise-constant interpolation operator, which includes in its range only one sample of the smooth error, produces efficiency degradation of our AMG method and more accurate interpolation operators obtained by weighted-Jacobi smoothing of the piecewise constant interpolation operators are required. This leads to a smoothed aggregation-type AMG method, which exhibits improved convergence and scalability properties, but it produces more dense coarse matrices and thus larger complexity, especially for linear systems arising from 3d PDE problems. Here we propose an extension of the method where more effective prolongator operators are defined, as presented next.
We follow the general strategy of the traditional unsmoothed aggregation methods when we want to incorporate several {\em algebraically smooth} vectors into the coarse hierarchy.
What is different in our approach is the origin of these vectors and the way we  build the hierarchy of aggregates. The details are outlined in the following section.

\section{Multiple-vector interpolation matrices in aggregation AMG}\label{section: multiple vectors interpolation}

To create aggregates of sufficiently large size we successively merge several levels of pairwise aggregates which we assume are already available.
In our case the hierarchy of such  pairwise aggregates are generated by the basic aggregation scheme discussed in the previous section.
The resulting sequence is denoted by
\begin{equation}
\A^k= \{ a_j^k \}_{j=1}^{n_k},\; k =1,2,\;\dots,nl.
\label{aggs}
\end{equation}
At each level k, the aggregates are represented by the following binary matrices:
\begin{equation}
\pi_{ij}^k= \left \{
\begin{array}{lr}
1 & \text{if } i \in a_j^k\\
0 & \text{otherwise}
\end{array}
\right. \ \ i=1, \ldots, n_{k-1}, \; j=1, \ldots, n_k,
\label{binmat}
\end{equation}
where $n_k$ is the number of aggregates at the level $k$ and $n_0=n$ is the number of original degrees of freedom (dofs). The above operator allows to map vectors associated with a coarse set of dofs ${1,2, \ldots, n_k}$ (the aggregates) on the finer set of dofs ${1,2, \ldots, n_{k-1}}$.
At the finest level, we have $\A^1= \{ a_j^1 \}_{j=1}^{n_1}$, where $a_j^1$ stands for aggregate of fine-level dofs and $n_1$ is the number of aggregates. The corresponding mapping operator is $\pi^1 \in \R^{n \times n_1}$.
The second assumption is that we have several fine-level samples of smooth vectors, $\{\bw_r\}^{n_w}_{r=0}$, that we will use to build
a coarse-level vector space defined as $\Range(P^1)$. These sample vectors can be computed by the application of the bootstrap procedure described in Algorithm~\ref{alg-boot}.

Let $k=1$ be and let $\left . \bw_r \right |_{a_j^1}$ be the $L_2$ orthonormal projection of $\bw_r$ on the aggregate $a_j^1$, we can write a rectangular matrix $P_{a_j^1} = \left [
\left . \bw_0 \right |_{a_j^1},\dots,\; \left .  \bw_{n_w} \right |_{a_j^1}\right ]$, which is of size $n_{a_j^1} \times (n_w+1)$, where $n_{a_j^1} > (n_w+1)$ is the size of the aggregate ${a_j^1}$. After forming $P_{a_j^1}$, we seek a local basis of the set of the given smooth vectors on the aggregate $a_j^1$ by performing a SVD of $P_{a_j^1}$ which gives:
\begin{equation}
\label{svd}
P_{a_j^1} = U_{a_j^1} \Sigma_{a_j^1} V^T_{a_j^1},
\end{equation}
where $U_{a_j^1}=[\bu_0, \ldots, \bu_{n_w}]$ has orthogonal columns and $\Sigma_{a_j^1}$ is the $(n_w+1) \times (n_w+1)$ diagonal matrix with the singular values on its main diagonal $\sigma_0 \geq \sigma_1 \ldots \geq \ldots \sigma_{n_w} \geq 0$.

After neglecting possible near-zero singular values of $\Sigma$, $\sigma_{n^{'}_{w}}, \ldots \sigma_{n_w} \leq TOL_{a_j^1}$, where $TOL_{a_j^1}$ is a given threshold for aggregate $a_j^1$, we use the first $n^{'}_{w}$ left singular vectors, i.e., the first $n^{'}_{w}$ columns of
$U_{a_j^1}$, to form a block $\tilde{P}_{a_j^1}$ of the following interpolation matrix:
\[
P^1=blockdiag(\tilde{P}_{a_1^1}, \ldots, \tilde{P}_{a_{n_1}^1}),
\]
where $\tilde{P}_{a_j^1}=[\bu_0, \ldots, \bu_{n^{'}_{w-1}}]$. Note that $P^1$ has orthogonal columns since each block $\tilde{P}_{a_j^1}$ has orthogonal columns, i.e., $\tilde{P}^T_{a_j^1}\tilde{P}_{a_j^1}=I$, and the same holds for $P^1$, i.e., $(P^1)^TP^1=I$. Therefore, by applying a standard Galerkin approach, we can compute a coarse-level matrix $A^1=(P^1)^TAP^1$. By construction, $P^1$ is a piecewise-constant interpolation operator, whose sparsity pattern is completely defined by the aggregates. The vector space $\Range(P^1)$ defines a coarse-level space which includes a set of smooth error samples of the given problem. We observe that the number $n^{'}_{w}$ of the left singular vectors in~\Ref{svd}, chosen as local basis of the smooth vectors for each aggregate, can be different for each aggregate, and the new coarse set of dofs has final dimension $n_c=\sum_{j=1}^{n_1} n^{'}_{w}$.

\subsection{Extension to the multi-level version}

We can iterate the above process for each new level $k=2, \ldots, nl$ of the hierarchy of aggregates in~\Ref{aggs}, generating a hierarchy of
block prolongators $P^k, \ k=1, \ldots, nl$ and respective coarse matrices $A^k$.

We observe that, for moving at a second coarse level it is needed to define a new binary transfer operator to map vectors from the first-level coarse set of dofs to the second-level set of dofs in~\Ref{aggs}, which is represented by the operator $\pi^2$ in~\Ref{binmat}.
To this aim we define the matrix $Q \in \R^{n_c \times n_1}$ with the following entries:
\begin{equation*}
Q_{ij}= \left \{
\begin{array}{lr}
1 & \text{if } i \in a_j^1\\
0 & \text{otherwise}
\end{array}
\right. \ \ i=1, \ldots, n_{c}, \; j=1, \ldots, n_1,
\end{equation*}
so that we can move at the second level of a new hierarchy by using the composite transfer operator $Q \pi^2 \in  \R^{n_c \times n_2}$. Similar transfer operators can be defined at each new level of the new hierarchy. 

For each new level $k = 2 \ldots, nl$, we consider as current-level smooth vectors the projection on the coarser level of the original $r$ smooth vectors, $\bw_r^k = (P^k)^T \ldots (P^1)^T \bw_r, \ $ $r = 0, \ldots, n_w$.

The setup of each new prolongator operator requires a number of $n_k$ small and dense SVD computations, where $n_k$ is the number of aggregates at each level. Therefore, using a very aggressive coarsening by merging an arbitrary large number of basic pairwise aggregates reduces the number of SVD computations at each level, although it increases the number of rows of each small matrix $P_{a_j^k}$.
However, it is worth noticing that SVD computations are completely independent on each aggregate, i.e., this part of the setup phase can be beneficial for parallel implementation.

The number of sample smooth vectors is another parameter of our method and increasing the number of sample smooth vectors generally improves the convergence (see Section~\ref{section: numerical results}). On the other hand, this increase requires that Algorithm~\ref{alg-boot} has to be applied for a larger number of iterations to estimate the needed smooth vectors and further it leads to larger SVD computations and to denser coarse matrices. Therefore, increasing the number of sample smooth vectors results in larger cost both for the setup and for the application of the final AMG method. We show that a good trade-off between the size of aggregates and the number of the smooth vectors, depending on the problem dimension, can lead to reliable AMG methods with reasonable operator complexity whose application in difficult problems is efficient and scalable.

Our procedure to compute the interpolation operator for multiple sample of smooth errors has some analogy with the method proposed in~\cite{C-06}. In that case, sample of smooth vectors are obtained during the AMG setup algorithm by applying a number of relaxation steps to the linear systems $A^{k-1}\bx = \b0$ for each new level $k = 1, \ldots nl-1$. Relaxation starts from a random vector at the finest level, while it starts from block vectors formed with blocks $\Sigma_{a_j^k} V^T_{a_j^k}$ obtained by the SVD in~\Ref{svd} for the aggregate $a_k^j$ at the successive levels. Our method differs from the method proposed in~\cite{C-06} also in the way it constructs the aggregates. In~\cite{C-06} the author applies the aggregation scheme proposed in~\cite{VMB-96}, using a measure of strength of connection among the variables, coupled with block prolongators based on interpolation of local basis for multiple vectors, with the final aim to represent more low-energy vectors further than the near-null space vector in the coarse space within a smoothed aggregation approach.

Here we adopt coarsening in which aggregates of arbitrary large dimension are obtained without referring to any standard measure of strength of connection among the variables well understood only for $M-$matrices. Dimensions of the aggregates are only related to the size of the original problem and the performance requirements. Indeed final computational complexity of the preconditioner is the result of a trade-off between size of aggregates and number of smooth vectors to be projected on them. No smoothing is applied to the final block prolongator to avoid further increase in the complexity of the final preconditioner. We refer to the thus modified  method as \emph{single-hierarchy multiple-vector aggregation-based AMG}. Observe that after the generation of the {\em algebraically smooth} vectors, we have available several hierarchies of pairwise aggregates (each one corresponding to one of the vectors), and we could use  any of them as starting set for the new aggregation to fit the  multiple vectors. In our experiments, we used the hierarchy computed at the last bootstrap stage which we found most appropriate.

\section{Numerical results}\label{section: numerical results}

In the present section we show results from implementing the \emph{single-hierarchy multiple-vector aggregation-based AMG} method discussed in this paper and compare it with our previously developed bootstrap AMG method, which is based on compatible matching as described in detail in the previous sections.
The comparison is done in terms of solve time, operator complexity (reflecting memory usage) and algorithmic scalability, all described in more details below.
The methods are tested as preconditioners in the Conjugate Gradient (CG) method.
For a given linear system, we set the unit vector as right-hand side and start with zero initial guess, Iterations are stopped when the euclidean norm of the residual is reduced by a factor of $10^{-6}$ or a maximum number of iterations $itmax=1000$ is reached.

The runs have been carried out on one core of a 2.6 GHz Intel Xeon E5-2670, running the Linux 2.6 kernel with the GNU compiler version
4.9. Release 5.0 of SuperLU~\cite{L-05}  is used for computing the LU factorization and triangular system solutions of the coarsest level systems of equations.
We use the Sparse Parallel Robust Algorithms Library (SPRAL)~\cite{HS-15}, implementing an auction-type algorithm, for computation of approximate
maximum product matchings. The latter is in the core of our \emph{coarsening based on compatible weighted matching}.

The algebraically smooth vectors $\{\bw_r\}$ are generated by the original bootstrap described in Algorithm~\ref{alg-boot}, starting from the initial
smooth vector $\bw_0 = {\mathbf 1}$. Note that we obtain very similar results also by starting from a random vector on which a sufficiently large number (e.g., 20) of symmetrized Gauss-Siedel relaxation steps is applied. The parameters of the bootstrap AMG algorithm are optimized to minimize the complexity of each cycle by balancing the coarsening factor achieved by merging two pairwise aggregation steps. For the actual details we refer to the implementation of our BootCMatch code described in detail in~\cite{DFV-18}. For the test cases discussed in this paper the new method composes three consecutive prolongators of the first hierarchy built by the original bootstrap in order to have sets of (sufficiently large) aggregates of size at most $64$.
All the AMG components, in the setup phase of the original bootstrap are applied as K-cycles (i.e., nonlinear AMLI cycle MG, see, \S~10.4 in~\cite{MLBFP}), with two inner iterations of flexible CG at each level but the coarsest one, in order to compensate for the unsmoothed aggregation and have good convergence components. The forward/backward Gauss-Seidel (GS) relaxation is used as pre-/post-smoother at all levels but the coarsest one, where a direct method is employed. The other characteristics of the bootstrap AMG are as follows:
\begin{itemize}
\item [(i)] we use the symmetrized multiplicative version, i.e, at stage $r$, we apply $\nu$-times  $(I-B_1A) \ldots (I-B_{r}A) (I-B_{r} A) \ldots (I-B_1A)$ to a random vector;
\item [(ii)] the number of steps at every stage was $\nu = 15$;
\end{itemize}
The above choices were experimentally found (in~\cite{DFV-18}) to lead to good trade-off between accuracy for computing the $r$th smooth vector $\bw_r$ and the computational cost.
For the computation of the blocks $\tilde{P}_{a_j^k}$ of the new multiple-vector prolongator after SVD computations of $P_{a_j^k}$, we use the left singular vectors such that the corresponding singular values are larger than $TOL_{a_j^k} = TOL * (size(a_j^k)/size(A))$, with $TOL=0.1$. The final \emph{single-hierarchy multiple-vector aggregation-based AMG} method is applied as a standard V-cycle to compensate the increasing operator complexity due to the use of multiple smooth vectors.

The following characteristics of the resulting single-hierarchy AMG method are reported in the Tables:
\begin{itemize}
\item number of levels \textbf{nl};
\item operator complexity:
\[
\textbf{opc}=\frac{\sum_{k=0}^{nl-1}nnz(A^{k})}{nnz(A^0)},
\]
where $A^k$ is the matrix at level $k$ ($k=0$ corresponds to the fine level) and $nnz(A^k)$ is the number of nonzeros of $A^k$; it gives information about the storage required by the operators of the AMG hierarchy as well as gives an estimate for the computational complexity in a standard V-cycle.
\item average coarsening factor defined as:
\[
\textbf{cr}= \frac{1}{nl} \sum_{k=1}^{nl} \frac{n(A^{k-1})}{n(A^{k})},
\]
where $n(A^k)$ is the size of $A^k$;
\item \textbf{tb}: the setup time in seconds, needed for building the AMG hierarchy;
\item \textbf{mvtb}: the part of \textbf{tb} needed for computing the hierarchy of block-diagonal prolongators $P^k$ and related coarse matrices, which also includes the time for creating large aggregates and the time for performing local SVD needed to extract linearly independent subset from the smooth vectors restricted to the current-level individual aggregates;
\item $\mathbf{\rho}$: an estimate of the convergence factor.
\end{itemize}
Number of iterations (\textbf{nit}) and execution time in seconds (\textbf{ts}) for the application of the preconditioned CG are also shown.
We did experiments for increasing number of samples of smooth vectors (\textbf{nsv}) from 3 till 10,
and generally report results when the corresponding single-hierarchy AMG shows solve time smaller than the original bootstrap AMG (see \textbf{b-it} and \textbf{b-ts} for number of iterations and solve times of the original bootstrap AMG). The original bootstrap AMG was applied in the default conditions, i.e., by using double pairwise unsmoothed aggregation coupled with K-cycle and symmetrized multiplicative composition of the AMG components. All the other algorithmic choices are the same as in the new method. The bootstrap process ends either when  the composite AMG reaches a convergence factor less than $\rho_{des}=0.8$ or when a maximum number of $15$ components are built.

\subsection{Anisotropic Diffusion}
\label{anisotropicPDE}

We started our analysis with test cases arising from the following anisotropic 2D PDE on the unit square coupled with homogeneous Dirichlet boundary conditions:
\[
- \div(K\; \nabla u)=f,
\]
where $K$ is the coefficient matrix
\[
K = \left [
\begin{array}{ll}
a & c\\
c & b
\end{array}
\right ], \quad
\text{ with } \quad \left \{
\begin{array}{l}
a= \epsilon + \cos^2(\theta)\\
b= \epsilon + \sin^2(\theta)\\
c=  \cos(\theta)\sin(\theta)
\end{array}
\right .
\]
The parameter $0 < \epsilon \leq 1$ defines the strength of anisotropy in the problem,  whereas the parameter $\theta$ specifies the direction of anisotropy. In the following we discuss results related to test cases with $\epsilon=0.001$ and $\theta= 0$, $\pi/8$, which we refer to as $ANI1$ and $ANI2$, respectively. The problem was discretized using the Matlab PDE toolbox, with linear finite elements on (unstructured) triangular meshes of three different sizes ($168577$, $673025$, $2689537$), obtained by uniform refinement. It is well known that the near-kernel of the above problems has dimension $1$ and includes the unit vector, then we expect that our single-hierarchy multiple-vector AMG shows good convergence already for a small number of additional smooth vectors. On the other hand, the original method performs very well also by using a single-hierarchy built on the base of the unit vector (see results in~\cite{DFV-18}, Table 3), indeed using bootstrap improves convergence rate at the cost of an increase both in setup and in solve time. We show here that using the new multi-vector single-hierarchy AMG allows further improvement in solve time with respect to the original AMG when no bootstrap is applied.

In Tables~\ref{ani-1}-\ref{ani-2}, we show results for increasing problem size and for increasing samples of smooth vectors.
We observe that adding few samples of smooth vectors to the starting unit one leads to an improvement in convergence rate and to a consequent reduction of iteration numbers and solve time, while for a number of smooth vectors larger than $6$ the convergence rate of the new method stabilizes. The best solve time is obtained when $5$ total vectors are used for the smallest and medium size problems and when $6$ vectors are used for the largest size, corresponding to AMG with a total operator complexity not larger than $2$. In about all cases reported in Table~\ref{ani-1}, solve time of the new method are also better than the original single-hierarchy AMG when no bootstrap is used and K-cycle is applied (see results in~\cite{DFV-18}, Table 3).
{\small
\begin{center}
\begin{table*}[t]
\caption{ANI1 test case for increasing size.\label{ani-1}}
\centering
\begin{tabular}{lccccc|cc}
& \multicolumn{5}{c}{\textbf{Setup}} & \multicolumn{2}{c}{\textbf{Solve}} \\
\cmidrule{2-6} \cmidrule{7-8}
\textbf{nsv} & \textbf{nl} & \textbf{opc} & \textbf{$\rho$} & \textbf{cr} &
\textbf{tb (mvtb)} & \textbf{nit}  & \textbf{ts} \\
\cmidrule{2-8}
& \multicolumn{5}{c}{\textbf{n=168577}} & \textbf{b-it=15} & \textbf{b-ts=4.25}\\\cmidrule{2-8}
3 & 3 & 1.46 & 0.91 & 7.49 & 11.95 (6.55) &  74  & 2.18\\
4 & 3 & 1.72 & 0.89 & 6.16 & 17.96 (9.51) &  50  & 1.77\\
5 & 3 & 2.04 & 0.88 & 5.29 & 25.26 (13.04)&  36  & 1.44\\
6 & 3 & 2.42 & 0.87 & 4.65 & 34.65 (17.88)&  34  & 1.74\\
7 & 3 & 2.85 & 0.86 & 4.20 & 45.67 (23.69)&  31  & 2.05\\
8 & 3 & 3.36 & 0.85 & 3.82 & 55.91 (28.24)&  28  & 1.83\\
9 & 3 & 3.89 & 0.84 & 3.56 & 69.64 (35.32)&  27  & 2.22\\
10& 3 & 3.40 & 0.85 & 3.81 & 80.7 (36.84) &  30  & 2.25\\
 \cmidrule{2-8}
& \multicolumn{5}{c}{\textbf{n=673025}} & \textbf{b-it=19} & \textbf{b-ts=28.22} \\\cmidrule{2-8}
3  &  3 & 1.38 & 0.91 & 8.47  & 78.58  (54.62) & 111 & 13.53 \\
4  &  3 & 1.60 & 0.90 & 7.13  & 118.45 (80.83) & 72  & 10.78 \\
5  &  3 & 1.86 & 0.89 & 6.22  & 167.86 (113.51)& 47  & 8.20\\
6  &  3 & 2.07 & 0.88 & 7.25  & 153.42 (79.82) & 51  & 8.32\\
7  &  3 & 2.54 & 0.87 & 5.13  & 329.71 (232.04)& 36  & 11.46\\
8  &  3 & 2.94 & 0.85 & 4.76  & 382.43 (259.12)& 34  & 11.06\\
9  &  3 & 3.41 & 0.85 & 4.47  & 472.30 (319.80)& 33  & 12.60\\
10 &  3 & 2.99 & 0.85 & 4.71  & 630.55 (351.56)& 37  & 14.46\\
\\\cmidrule{2-8}
& \multicolumn{5}{c}{\textbf{n=2689537}} & \textbf{b-it=22} & \textbf{b-ts=168.21}\\\cmidrule{2-8}
3  &  3 & 1.32 & 0.91 & 13.65 & 341.23  (235.46) & 202 & 101.30 \\
4  &  3 & 1.51 & 0.90 & 12.17 & 550.90  (382.48) & 197 & 140.94 \\
5  &  3 & 1.73 & 0.89 & 11.25 & 764.42  (522.16) & 141 & 112.08\\
6  &  3 & 2.00 & 0.88 & 10.63 & 989.91  (662.49) & 81  & 57.14\\
7  &  3 & 2.31 & 0.88 & 10.17 & 1452.66 (987.20) & 73  & 73.75\\
8  &  3 & 2.66 & 0.87 & 9.82  & 1696.49 (1138.66)& 68  & 78.29\\
9  &  3 & 2.32 & 0.88 & 9.94  & 3570.05 (1122.38)& 67  & 67.69\\
10 &  3 & 2.68 & 0.87 & 9.63  & 3998.39 (1460.57)& 68  & 79.50\\
\end{tabular}
\end{table*}
\end{center}
}

\begin{figure*}[htb]
\begin{center}
\includegraphics[width=0.45\textwidth]{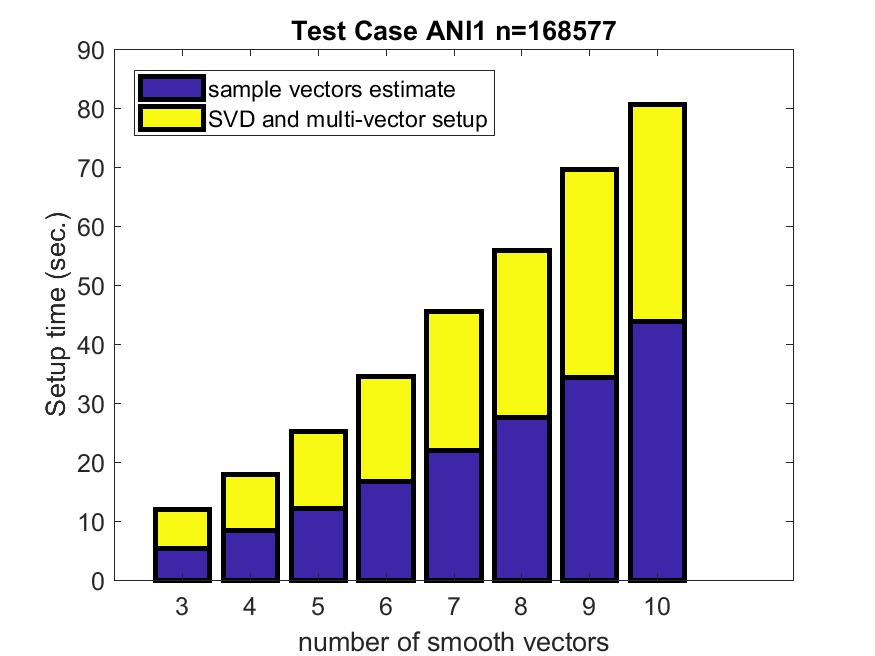}
\includegraphics[width=0.45\textwidth]{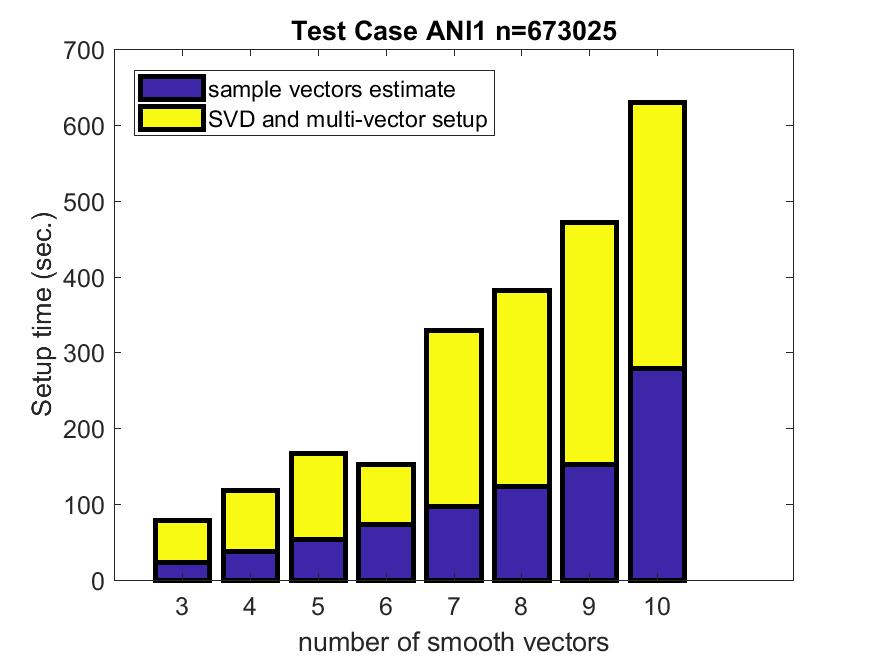}
\includegraphics[width=0.45\textwidth]{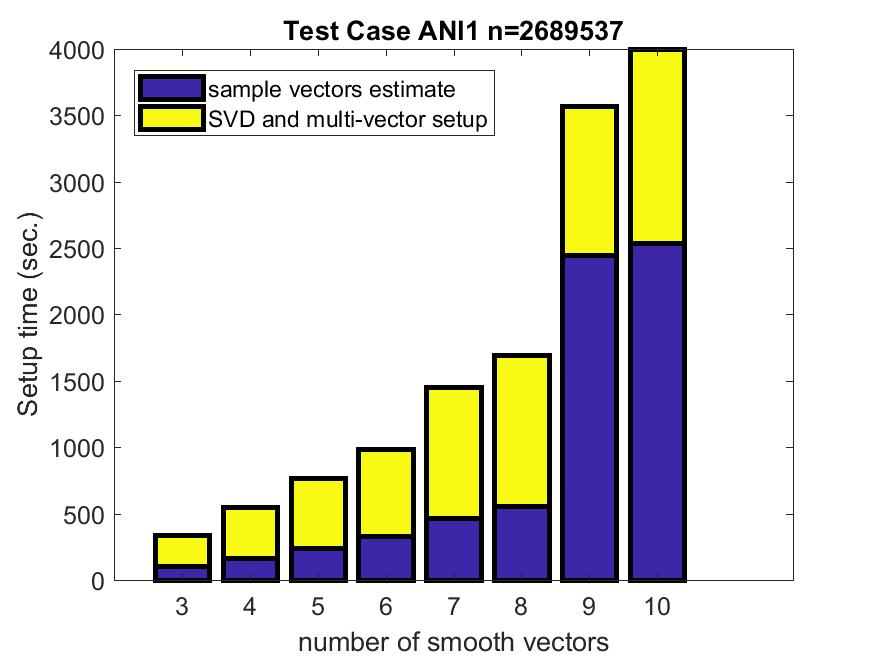}
\end{center}
\caption{ANI1: Setup time of the single-hierarchy multiple-vector AMG\label{fig-ANI1}}
\end{figure*}

As expected, increasing number of  smooth vectors leads to an increase of operator complexity with a corresponding decrease in the coarsening ratio, due to prolongators having more columns and then producing larger coarse matrices.
Furthermore, setup time of the new method increases for increasing number of smooth vectors. It includes the time to generate smooth vectors by bootstrap and the time for SVD computations needed to build linear independent projections of the smooth vectors on the aggregates.
We observe that for a fixed number of smooth vectors the percentage of SVD computations with respect to the total setup times increases for increasing matrix sizes, indeed for aggregates of a fixed size (at most $64$, in our experiments), the number of aggregates at each level increases with matrix size. On the other hand, when the number of smooth vectors increases, the percentage of SVD computations on the total setup time decreases.
This behavior is clearly illustrated for the ANI1 test case in Fig.~\ref{fig-ANI1}, where we describe the total setup time (\textbf{tb}) defined as the sum of the time spent in the bootstrap process needed for generating the smooth vectors, and the time spent in the SVD and in building the new multi-vector prolongators (\textbf{mvtb}).
It is worth noting  that SVD and the setup of the multi-vector prolongators are local (aggregate-by-aggregate) and hence embarrassingly parallel, and that cost can be amortized for multiple r.h.s.
We finally note that the original bootstrap AMG needs to generate $6$, $7$ and $8$ hierarchies to reach the desired convergence rate, requiring setup times equal to $16.70$, $99.72$ and $569.66$ for the three problem sizes reported in Table~\ref{ani-1}, respectively. Therefore, in the cases of using $3$ and $4$ smooth vectors, the new method shows total times which are smaller than the total time of the built bootstrap AMG, with significant memory savings.
More specifically, for the largest size matrix, the operator complexity of the hierarchy obtained with $4$ vectors is $1.51$, while the original bootstrap AMG shows average operator complexity of about $1.4$ for all hierarchies; therefore, in the case of $4$ hierarchies built employing the first 4 vectors, the original bootstrap requires a memory increase of a factor about $5.6$ with respect to the memory needed for the system matrix and shows a convergence rate of about $0.91$, corresponding to solve time of $250.94$ and $57$ iterations. In conclusion, using the multiple-vector prolongators allows us to obtain a single-hierarchy that, applied as a simple V-cycle, has convergence rate, memory requirements and total time better than the composite bootstrap AMG with the same number of smooth vectors and shows a moderate increase of iteration number when the problem size increases. On the other hand, further reduction of solve time can be obtained at the expense of an increase in setup times and a moderate increase in memory requirements.
{\small
\begin{center}
\begin{table*}[t]
\caption{ANI2 test case for increasing size.\label{ani-2}}
\centering
\begin{tabular}{lccccc|cc}
& \multicolumn{5}{c}{\textbf{Setup}} & \multicolumn{2}{c}{\textbf{Solve}} \\
\cmidrule{2-6} \cmidrule{7-8}
\textbf{nsv} & \textbf{nl} & \textbf{opc} & \textbf{$\rho$} & \textbf{cr} &
\textbf{tb (mvtb)} & \textbf{nit}  & \textbf{ts} \\
\cmidrule{2-8}
& \multicolumn{5}{c}{\textbf{n=168577}} & \textbf{b-it=16} & \textbf{b-ts=4.49}\\\cmidrule{2-8}
3 & 3 & 1.46 & 0.90 & 7.52 & 11.84 (6.50) &  86  & 2.46\\
4 & 3 & 1.71 & 0.89 & 6.20 & 17.74 (9.33) &  50  & 1.70\\
5 & 3 & 2.03 & 0.88 & 5.32 & 24.98 (12.86)&  42  & 1.58\\
6 & 3 & 2.41 & 0.87 & 4.6  & 34.32 (17.69)&  40  & 1.88\\
7 & 3 & 2.84 & 0.86 & 4.21 & 44.38 (22.62)&  32  & 1.81\\
8 & 3 & 3.33 & 0.84 & 3.84 & 56.57 (28.95)&  28  & 1.96\\
9 & 3 & 2.91 & 0.85 & 4.15 & 67.65 (30.51)&  33  & 2.28\\
10& 3 & 3.41 & 0.84 & 3.79 & 80.28 (36.49)&  26  & 1.78\\
 \cmidrule{2-8}
& \multicolumn{5}{c}{\textbf{n=673025}} & \textbf{b-it=15} & \textbf{b-ts=25.14} \\\cmidrule{2-8}
3  &  3 & 1.38 & 0.91 & 8.49  & 78.36  (54.60) & 124 & 15.01 \\
4  &  3 & 1.60 & 0.90 & 7.10  & 116.63 (79.41) & 70  & 9.58 \\
5  &  3 & 1.86 & 0.89 & 6.22  & 165.00 (111.37)& 58  & 9.26\\
6  &  3 & 2.17 & 0.87 & 5.60  & 225.68 (152.23)& 42  & 8.31\\
7  &  3 & 2.54 & 0.86 & 5.11  & 296.39 (200.01)& 40  & 10.04\\
8  &  3 & 2.94 & 0.86 & 4.76  & 376.26 (254.11)& 42  & 12.97\\
9  &  3 & 2.56 & 0.87 & 5.08  & 503.14 (258.15)& 46  & 15.12\\
10 &  3 & 2.99 & 0.86 & 4.71  & 581.16 (307.72)& 42  & 11.69\\
\\\cmidrule{2-8}
& \multicolumn{5}{c}{\textbf{n=2689537}} & \textbf{b-it=19} & \textbf{b-ts=237.85}\\\cmidrule{2-8}
3  &  3 & 1.32 & 0.91 & 13.61 & 341.69  (236.76) & 214 & 107.48 \\
4  &  3 & 1.50 & 0.90 & 12.16 & 512.53  (347.35) & 193 & 104.86 \\
5  &  3 & 1.73 & 0.90 & 11.33 & 717.12  (479.33) & 137 & 81.51\\
6  &  3 & 1.99 & 0.89 & 10.70 & 1006.93 (663.61) & 96  & 78.21\\
7  &  3 & 2.30 & 0.88 & 10.14 & 1280.30 (853.87) & 74  & 57.49\\
8  &  3 & 2.66 & 0.88 & 9.80  & 1654.71 (1073.02)& 65  & 61.08\\
9  &  3 & 2.32 & 0.88 & 9.97  & 3640.83 (1229.98)& 79  & 78.90\\
10 &  3 & 2.68 & 0.88 & 9.62  & 4063.48 (1522.22)& 77  & 87.18\\
\end{tabular}
\end{table*}
\end{center}
}
In Table~\ref{ani-2} we have similar results for the more challenging ANI2 test case, with anisotropy direction not aligned with Cartesian axes. In this case, the best solve time is obtained when $5$, $6$ and $7$ smooth vectors are employed for the three matrix size, respectively. On the other hand, the original bootstrap AMG requires
$6$, $8$ and $9$ hierarchies to reach the desired convergence rate, with setup times equal to $16.59$, $125.17$ and $843.34$ for the three problem sizes, respectively. Therefore, for the largest matrix size, the new method is able to obtain a total time smaller than the original bootstrap AMG by using 3 to 6 smooth vectors with a largely significant memory savings. Indeed, in the case of $6$ hierarchies, the bootstrap AMG which also in this test case shows average operator complexity of about $1.4$ for all hierarchies, requires a memory increase by a factor of about $8.4$ with respect to the memory needed for the system matrix, showing a convergence rate of about $0.87$ and solve time of $294.53$ for $37$ iterations, while the new multiple-vector AMG with 6 smooth vectors shows an operator complexity less than $2$ which is a significant improvement.

\subsection{Linear elasticity}
\label{linearelasticityPDE}

\begin{figure}[htb]
\begin{center}
\includegraphics[width=0.4\textwidth]{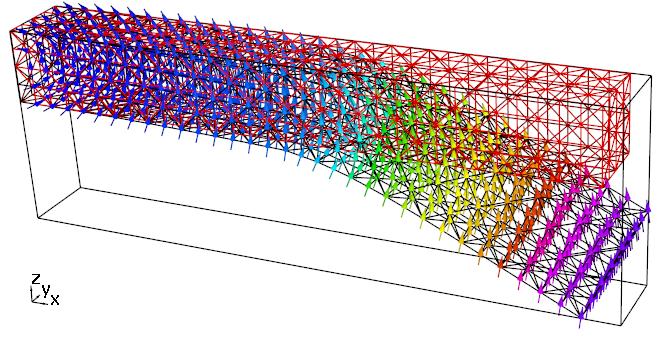}
\includegraphics[width=0.4\textwidth]{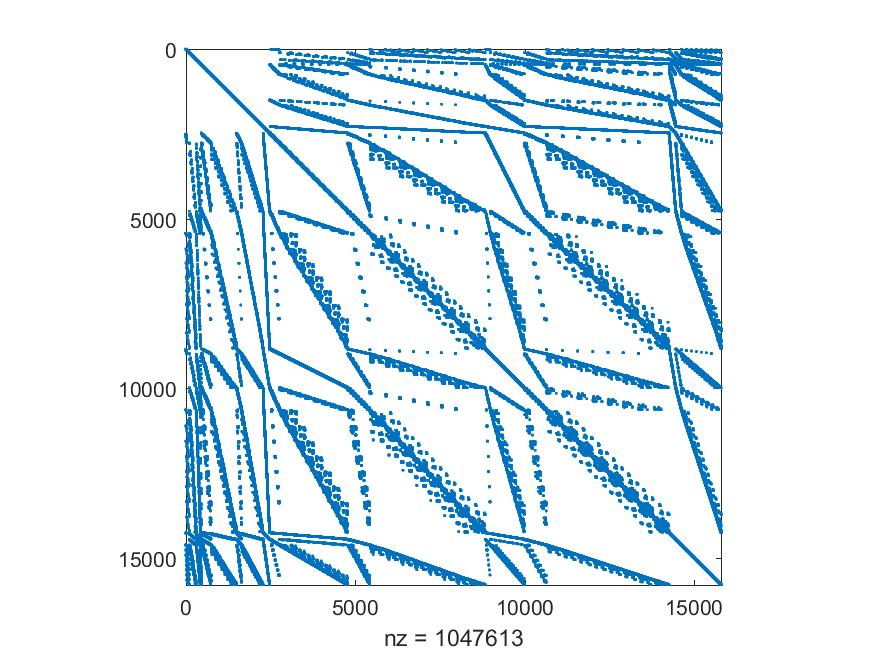}
\end{center}
\caption{Linear elasticity test case: mesh geometry (left); sparsity pattern of the system matrix (right)~\label{beam}}
\end{figure}

In this section we discuss results for test problems arising from discretization of the Lam\'e equations of linear elasticity:
\[
\mu \Delta \bu +(\lambda + \mu) \nabla(\div \bu)= \mathbf{f} \quad \quad  \bx \in \Omega \subset {\mathbb R}^d,
\]
where $\bu=\bu(\bx)$ is the displacement vector, $\Omega$ is the
3d spatial domain illustrated in Fig.~\ref{beam}, and $\lambda$ and $\mu$ are the Lam\'e constants. It is well-known that, when $\lambda >> \mu$, i.e., when the material is nearly incompressible, the problem becomes very ill-conditioned. We consider Lam\'e equations on a long beam having an aspect ratio of 1:8,  characterized by $\mu=0.5$ and two different values of
$\lambda_1=7, \; \lambda_2=10$, corresponding to two test cases of increasing conditioning; one side of the beam is considered fixed and the opposite
end is pushed downward, i.e., mixed Dirichlet and traction-type conditions are applied (see Fig.~\ref{beam}).
The problem was discretized with linear finite elements on tetrahedral meshes, using the software package MFEM~\cite{mfem}; different problem sizes were obtained by uniform refinements. More specifically, for each test case, three different sizes are considered ($15795$,$111843$, $839619$).
We applied a {\em node-based}~\cite{R-86} discretization of the above PDE, where at each mesh point all scalar components of the displacement vector are grouped together. The above ordering leads to systems of equations whose coefficient matrix is s.p.d. and have small block entries of dimension $3 \times 3$. We referred to the system matrices as $LE1$ and $LE2$, corresponding to the two increasing values of $\lambda$, respectively.

We observe that AMG methods, when applied to vector equations, employ ordering-aware coarsening schemes and relaxation methods to be effective~\cite{RS-87,BKY-10}. Furthermore, the construction of efficient AMG for linear elasticity relies on a priori knowledge of the rigid body modes ($6$ vectors in 3d problems), and requires that they are well represented on all coarse levels (see e.g., ~\cite{BKY-10}). Some other efficient AMG methods for finite element discretization of linear elasticity problems, namely AMG using element interpolation (AMGe) and its variants, require access to the fine-grid element matrices~\cite{JV-01,Cetal-03,KV-06,K-08}.
As we have demonstrated in a previous publication (\cite{DFV-18}) our bootstrap AMG is able to successfully handle the elasticity equations without using any a priori information on the rigid body modes and with no information on the discretization mesh and/or ordering. In what follows, as we can see from  Tables~\ref{mfem-1}-\ref{mfem-2}, our modified single-hierarchy AMG does as well. Additionally, it shows superior performance than the original multiple-component bootstrap AMG, which will become evident after a closer look at the tables.

In the case of LE1, the original bootstrap AMG generates $9$, $11$ and $14$ smooth vectors and setup times equal to $14.94$, $190.13$ and $2900.51$, for the three problem sizes, respectively.
As expected, for these vector linear systems arising from linear elasticity, a larger number of smooth vectors are needed to obtain
a multiple-vector AMG with good convergence and smaller solve time than the original bootstrap. On the other hand, increasing smooth vectors generally improves convergence rate and for the smallest size matrix, the new method obtains a convergence rate less than the original bootstrap already employing $8$ vectors. The best solve time are obtained by using $9$ vectors for the small size matrix, while using $10$ smooth vectors for the medium and largest matrix improves both convergence rate and solve time.

\begin{figure*}[htb]
\begin{center}
\includegraphics[width=0.45\textwidth]{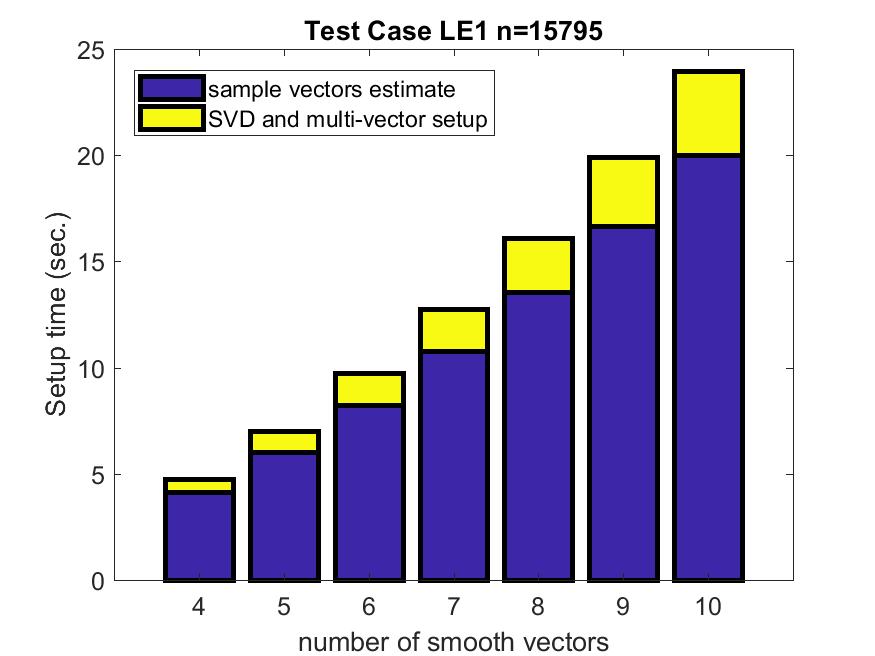}
\includegraphics[width=0.45\textwidth]{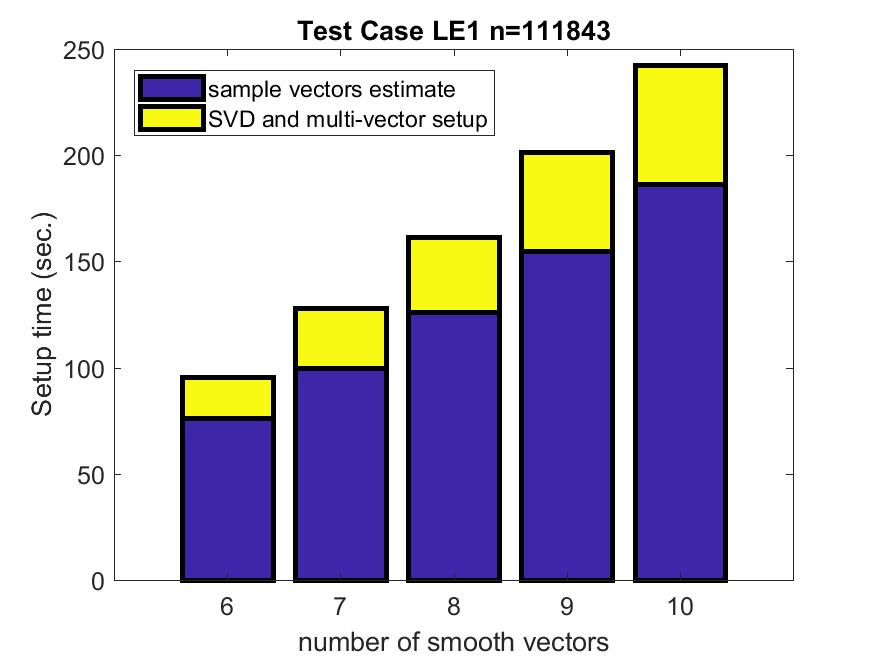}
\includegraphics[width=0.45\textwidth]{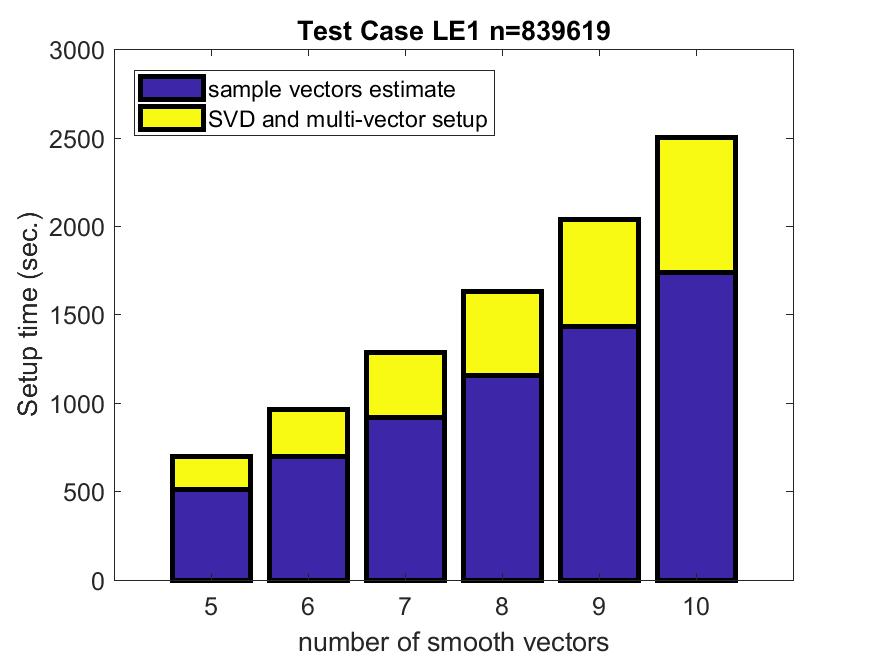}
\end{center}
\caption{LE1: Setup time of the single-hierarchy multiple-vector AMG\label{fig-LE1}}
\end{figure*}

Also in this case, for all problem sizes we can obtain a single-hierarchy AMG with smaller total times than the original bootstrap AMG by a suitable choice of the number of smooth vectors. On the other hand, if we look at the problem with the largest size, the new method, employing $10$ smooth vectors is already able to obtain total time much smaller than the original bootstrap AMG with a significant reduction in memory requirements. Indeed, the bootstrap AMG with the desired convergence rate ($0.8$) has an average operator complexity of about $1.4$ for all hierarchies, then with its $14$ built hierarchies requires a memory increase of a factor about $20$ with respect to the memory needed for the system matrix, while the new multiple-vector AMG with $10$ smooth vectors shows an operator complexity less than $5.2$. On the other hand, the original bootstrap using only $10$ hierarchies (and corresponding $10$ smooth vectors) shows a convergence rate of about $0.96$ and solve time of $2066.31$ for $102$ iterations.
{\small
\begin{center}
\begin{table*}[t]
\caption{LE1 for increasing size.\label{mfem-1}}
\centering
\begin{tabular}{lccccc|cc}
& \multicolumn{5}{c}{\textbf{Setup}} & \multicolumn{2}{c}{\textbf{Solve}} \\
\cmidrule{2-6} \cmidrule{7-8}
\textbf{nsv} & \textbf{nl} & \textbf{opc} & \textbf{$\rho$} & \textbf{cr} &
\textbf{tb (mvtb)} & \textbf{nit}  & \textbf{ts} \\\cmidrule{2-8}
& \multicolumn{5}{c}{\textbf{n=15795}} & \textbf{b-it=14} & \textbf{b-ts=2.01}\\
\cmidrule{2-8}
4  &  2 & 1.62 & 0.90 & 6.92 &  4.77  (0.62) & 96 & 1.67\\
5  &  2 & 2.03 & 0.87 & 5.63 &  7.02  (0.99) & 62 & 1.33\\
6  &  2 & 2.49 & 0.86 & 4.62 &  9.75  (1.49) & 42 & 1.14\\
7  &  2 & 3.07 & 0.81 & 4.03 &  12.76 (1.96) & 28 & 0.83\\
8  &  2 & 3.45 & 0.77 & 3.49 &  16.09 (2.52) & 23 & 0.82\\
9  &  2 & 3.98 & 0.73 & 3.18 &  19.93 (3.26) & 18 & 0.75\\
10  &  2 & 4.51 & 0.73 & 2.99 &  23.96 (3.97) & 16 & 0.75\\
\cmidrule{2-8}
& \multicolumn{5}{c}{\textbf{n=111843}} & \textbf{b-it=12} & \textbf{b-ts=20.33} \\
\cmidrule{2-8}
6  & 3 & 2.49 & 0.88 & 2.49 &  95.65 (19.52)& 92 & 19.46\\
7  & 3 & 3.07 & 0.87 & 3.07 & 128.05 (28.18)& 68 & 17.90\\
8  & 3 & 3.70 & 0.85 & 3.70 & 161.45 (35.42)& 47 & 14.35\\
9  & 3 & 4.41 & 0.83 & 4.41 & 201.50 (46.61)& 36 & 13.22\\
10  & 3 & 4.83 & 0.81 & 4.83 & 242.37 (56.07)& 27 & 11.30\\
\cmidrule{2-8}
& \multicolumn{5}{c}{\textbf{n=839619}} & \textbf{b-it=14} & \textbf{b-ts=314.50}\\\cmidrule{2-8}
5  & 3 & 2.00 & 0.89 & 8.91 & 701.71  (187.26) & 206 & 301.70\\
6  & 3 & 2.48 & 0.88 & 7.98 & 965.97  (264.77) & 148 & 260.90\\
7  & 3 & 3.04 & 0.87 & 7.38 & 1286.28 (367.85) & 106 & 226.00\\
8  & 3 & 3.64 & 0.86 & 6.85 & 1631.67 (472.07) & 76  & 189.27\\
9  & 3 & 4.38 & 0.84 & 6.48 & 2038.73 (604.47) & 59  & 173.53\\
10  & 3 & 5.12 & 0.83 & 6.20 & 2503.55 (761.74) & 46  & 158.93\\
\end{tabular}
\end{table*}
\end{center}
}
In Fig.~\ref{fig-LE1} we again describe the total setup time (\textbf{tb}) showing the percentage spent in generating the smooth vectors and that spent in SVD and in building the new multi-vector prolongators (\textbf{mvtb}). The general behavior is the same as in the ANI test cases. On the other hand, as expected, the time spent in the bootstrap procedure needed for generating the smooth vectors has a larger impact on the overall setup time for the LE test cases. Indeed in 3d PDE problems we have denser matrices at each level of the hierarchy which results in a more expensive application of each step of the bootstrap process, whereas we observe that using the new single-hierarchy multi-vector AMG allows to obtain better solve time than the original bootstrap AMG already using a small number of smooth vectors.
Finally, in the case of LE2 case, the original bootstrap AMG generates $11$, $12$ and $12$ smooth vectors and setup times equal to $19.92$, $222.05$ and $2252.99$, for the three problem sizes, respectively. The general behavior is similar to the previous test cases, confirming a very significant reduction of memory requirements and solve times of the new method with respect to the original bootstrap. In this more complex case we can observe that for increasing matrix size, the gain in solve time of the new method with respect to the original bootstrap appears also more significant for increasing number of smooth vectors.

We point out that we test some popular  AMG methods on the above linear elasticity test cases and our experiments show that AGMG~\cite{N-10}
has a divergent behavior, while BoomerAMG~\cite{HY-02} needs unknown-based discretization for convergence, while our approach is completely algebraic and independent of a priori information on the problem and of the discretization.
{\small
\begin{center}
\begin{table*}[t]
\caption{LE2 for increasing size.\label{mfem-2}}
\centering
\begin{tabular}{lccccc|cc}
& \multicolumn{5}{c}{\textbf{Setup}} & \multicolumn{2}{c}{\textbf{Solve}} \\
\cmidrule{2-6} \cmidrule{7-8}
\textbf{nsv} & \textbf{nlev} & \textbf{opc} & \textbf{$\rho$} & \textbf{cr} &
\textbf{tb (mvtb)} & \textbf{nit}  & \textbf{ts} \\\cmidrule{2-8}
& \multicolumn{5}{c}{\textbf{n=15795}} & \textbf{b-it=11} & \textbf{b-ts=1.95}\\\cmidrule{2-8}
5 & 2 & 2.04 & 0.89 & 5.64 &  7.07  (1.00) & 85 & 1.82\\
6 & 2 & 2.47 & 0.86 & 4.74 &  9.55  (1.33) & 58 & 1.42\\
7 & 2 & 3.01 & 0.84 & 4.05 &  12.67 (1.96) & 40 & 1.20\\
8 & 2 & 3.65 & 0.85 & 3.56 &  15.97 (2.52) & 30 & 1.05\\
9 & 2 & 4.37 & 0.84 & 3.17 &  20.20 (3.59) & 24 & 1.04\\
10 & 2 & 4.94 & 0.76 & 2.94 &  24.37 (4.36) & 22 & 1.07\\
\cmidrule{2-8}
& \multicolumn{5}{c}{\textbf{n=111843}} & \textbf{b-it=11} & \textbf{b-ts=20.36} \\\cmidrule{2-8}
9 & 3 & 4.61 & 0.86 & 3.95 & 205.79 (50.27) & 48 & 18.76 \\
10 & 3 & 5.89 & 0.85 & 3.68 & 258.84 (70.75) & 39 & 19.03\\
\cmidrule{2-8}
& \multicolumn{5}{c}{\textbf{n=839619}} & \textbf{b-it=45} & \textbf{b-ts=910.09}\\\cmidrule{2-8}
4   & 3 & 1.70 & 0.91 & 10.41& 489.30 (131.99) & 432 & 594.90\\
5   & 3 & 2.18 & 0.90 & 9.00 & 725.84 (205.60) & 283 & 450.58\\
6   & 3 & 2.65 & 0.90 & 8.02 & 995.69 (284.81) & 177 & 335.45\\
7   & 3 & 3.22 & 0.89 & 7.39 & 1323.90 (393.11)& 125 & 283.76\\
8   & 3 & 3.84 & 0.88 & 6.91 & 1683.30 (506.05)& 83  & 222.01\\
9   & 3 & 4.65 & 0.87 & 6.55 & 2111.70 (657.24)& 64  & 203.13\\
10  & 3 & 5.64 & 0.86 & 6.17 & 2613.75 (849.79)& 51  & 192.23\\
\end{tabular}
\end{table*}
\end{center}
}

\section{Conclusions}\label{section: conclusions}
In this paper we proposed a modification of a bootstrap AMG aimed to improve the computational costs in the solution phase of the resulting AMG. This is achieved by an aggressive coarsening algorithm which incorporates multiple algebraically smooth vectors generated by the original bootstrap algorithm, in a single hierarchy employing sufficiently large aggregates. These aggregates are in turn compositions of aggregates already built throughout the original bootstrap algorithm.
We have shown that the new \emph{single-hierarchy multiple-vector aggregation-based AMG} method has an overall better efficiency both in terms of memory and solve times. This is demonstrated on a class of linear systems arising from discretization of scalar and vector PDEs. The capability of the new method of defining aggregates of arbitrary size and by choosing an arbitrary number of samples of the created smooth vectors for the multiple-vector interpolation operators provides flexibility to balance the trade-off between computational complexity and convergence properties for the resulting AMG, depending on the specific application and computer platform features. Work  in progress involves the development of efficient parallel implementation of the code on multicore computers, including also GPU accelerators.
Furthermore, future work will focus on reducing setup time needed to generate accurate estimation of sample of smooth vectors by employing additive composition of the multiple components generated by the  bootstrap in order to enhance parallelism also in this phase.

\section{Acknowledgments}
\label{ack}

The authors wish to thank the anonymous reviewers for their comments which helped us to obtain an improved version of the paper.

\end{document}